\documentclass[10pt]{article}

\usepackage{amsthm}
\usepackage{amssymb}
\usepackage{amsmath}
\usepackage{enumerate}
\usepackage[mathscr]{eucal}

\ifx\pdfoutput\undefined
   \usepackage[dvips]{graphicx}
   \else
   \usepackage[pdftex]{graphicx}
   \fi
   \usepackage{epstopdf}

\theoremstyle{plain}
\newtheorem{tetel}{Theorem}[section]
\newtheorem{all}[tetel]{Proposition}

\newtheorem{kov}[tetel]{Corollary}

\theoremstyle{definition}\newtheorem{Def}[tetel]{Definition}
\theoremstyle{remark}\newtheorem{megj}[tetel]{Remark}
\newtheorem{pelda}[tetel]{Example}

\newcommand*{\Z}{\ensuremath{\mathbf Z}}
\newcommand*{\R}{\ensuremath{\mathbf R}}
\newcommand*{\di}{\mathrm d}

\begin{document}

\DeclareGraphicsExtensions{.eps,.pdf}

\title{Rulings of Legendrian knots as spanning surfaces}
\author{Tam\'as K\'alm\'an\\ University of Southern California / University of Tokyo}
\maketitle

\begin{abstract}
Each ruling of a Legendrian link can be naturally treated as a
surface. For knots, the ruling is $2$--graded if and only if the surface is orientable. For $2$--graded rulings of homogeneous (in particular, alternating) knots,
we prove that the genus of this surface is at most the genus of the knot. While this is not true in general, we do prove that the canonical genus (a.k.a.\ diagram genus) of any knot is an upper bound for the genera of its $2$--graded rulings.
\end{abstract}


\section{\small Introduction} 

A \emph{Legendrian knot} in $\R^3_{xyz}$ is a smooth
embedding of $S^1$ that is always tangent to the $2$--planes of
the so-called \emph{standard contact structure} $\xi=\ker(\di z-y\di x)$, and a \emph{link} is a finite disjoint union of knots.
We consider Legendrian knots up to \emph{Legendrian isotopy}, that is
homotopy through Legendrian knots. The most basic, ``classical''
invariants of Legendrian isotopy are the Thurston--Bennequin number
$tb$, rotation number $r$, and smooth type.\footnote{For definitions,
see \cite{etn}, or section 2 of \cite{en}. In general, we assume
that the reader is familiar with the basics of Legendrian knots
such as, for example, front diagrams.} (Of course, any invariant of
smooth isotopy is also an invariant of Legendrian isotopy.) 
In this paper, we will study the so-called ruling invariants, introduced by Chekanov and Pushkar \cite{chp} and independently Fuchs
\cite{fuchs}.

The beginning of Legendrian knot theory dates back to 1984 when Bennequin \cite{benn} proved the famous inequality
\begin{equation}\label{benn}tb(L)+|r(L)|\le2g(L)-1,\end{equation}
where $L$ is an arbitrary Legendrian knot and $g(L)$ is its (smooth) Seifert genus.
Our goal is to strengthen the relationship between the genus and Legendrian invariants by proving the following:

\begin{tetel}\label{thm:canon}
The genus of any $2$--graded ruling of the Legendrian
knot $L$ is less than or equal to the canonical genus of $L$.
\end{tetel}

Rulings and the idea that a genus may be associated to them are central to this paper. However the discussion of these notions is deferred to section \ref{sec:rul}. Recall that the \emph{canonical genus} $\tilde g$ is the minimum of the genera of spanning surfaces obtained using Seifert's algorithm on diagrams of the given knot. The proof of Theorem \ref{thm:canon}, given toward the end of section \ref{sec:dan}, is a trivial application of
a result of Rutherford \cite{rulpoly} and Morton's inequality \cite{morton}
that the $z$--degree of the Homfly polynomial bounds the canonical genus from below.

Because for homogeneous knots 
(a class of knots that includes all alternating and positive knots) the canonical genus agrees with the genus \cite[section 7.6]{crom}, the following is immediate.

\begin{kov}\label{kov}
The genus of any $2$--graded ruling of the homogeneous Legendrian
knot $L$ is less than or equal to the genus of $L$.
\end{kov}

In addition to homogeneous knots, the statement of the corollary holds for all prime knots up to $13$ crossings; 
see section \ref{sec:jim}. 
It also holds for connected sums of the aforementioned by the additivity of genus and Proposition \ref{pro:connsum}. It does not, however, hold in general. 

%
%

\begin{pelda}\label{ellenpelda}
There exist (prime) knots whose genus is less than the maximum of the genera of their $2$--graded rulings. Two examples are shown in Figure \ref{fig:14n22180}. One of these ($14n_{22180}$ in the Hoste--Thistlethwaite table) is a Whitehead double of the trefoil, but the other ($14n^*_{19265}$) is not a satellite knot (thus it is hyperbolic). 
Both are of genus $1$ and both possess rulings with genera up to $2$. The genus $2$ ruling is unique and $\Z$--graded for both. This shows that Corollary \ref{kov} fails in general for non-homogeneous knots even if we restrict the statement to $\Z$--graded rulings. 
\begin{figure}
   \centering
   \includegraphics[width=\linewidth,height=2.5in]{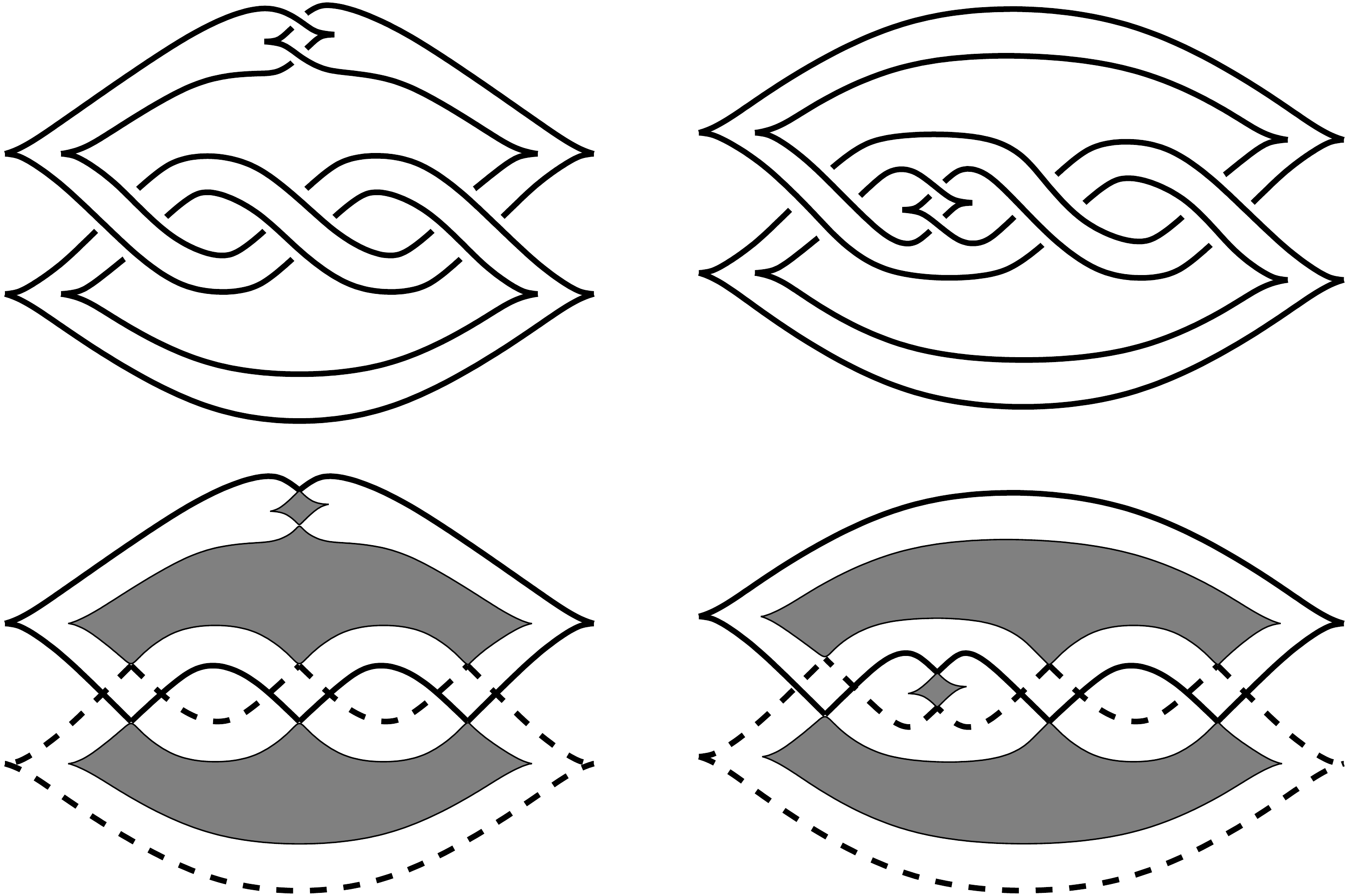}
   \caption{Knots with rulings whose genera exceed the knot's smooth genus. (The diagrams are non-generic but no switch is vertically aligned with a non-switching crossing.)}
   \label{fig:14n22180}
\end{figure}
\end{pelda}

\begin{megj}
A strengthening of our statements to the $4$--ball genus is not conceivable either. There exist slice knots, for example the alternating knot $10_{22}$ in Rolfsen's table, with positive genus $2$--graded rulings. 
\end{megj}

The paper is organized as follows. In section \ref{sec:rul}, we review rulings, their gradedness and genera, and establish some relations of these to the topology of the ruling. In section \ref{sec:dan}, we review Rutherford's work and prove Theorem \ref{thm:canon}. Section \ref{sec:jim} details how the two knots of Example \ref{ellenpelda} were found.

Acknowledgements: The phenomena described in this paper are related to Lagrangian surfaces in $B^4$ that a Legendrian in $S^3$ bounds. Thus this paper can be thought of as a prelude of forthcoming joint work \cite{cob} with Tobias Ekholm and Ko Honda. I had many inspiring discussions with them. I also thank Jim Hoste for the tremendous insight gained from the computer searches of knot tables that he conducted for me. My conversations with Francis Bonahon and Dan Rutherford were very useful, too.

\section{\small Rulings and their genera}\label{sec:rul}

The \emph{front diagram} of a Legendrian link $L\subset\R^3_{xyz}$ is its projection onto the $xz$--plane. For generic $L$, it has an equal number of so-called \emph{left and right cusps} (with respect to the $x$--direction), and the rest of the diagram consists of strands that connect a left cusp to a right cusp. A \emph{Maslov potential} \cite{chp,fuchs} is a locally constant $\Z_{2r}$--valued\footnote{The \emph{rotation} $r$, up to sign, of a Legendrian link is defined as the greatest common divisor of the rotations of its constituent components.} function that is defined along these strands with the following property: At each cusp, out of the two strands that join there, the one that is (locally) upper has potential that is $1$ higher than that of the other strand.

Recall that the \emph{index} of a crossing in a front diagram is the difference (in $\Z_{2r}$) of the Maslov potentials of the two intersecting strands. (The potential of the lower strand is to be subtracted from that of the upper one.) We adopt the convention 
that when a front diagram is oriented, strands traveling to the right (increasing $x$) have even Maslov potential. This implies that even for multi-component links, a crossing has even index if and only if it is positive. For the following definition \cite{chp,fuchs}, we assume that the front diagram of the Legendrian is generic in that no pair of crossings share the same $x$--coordinate. 

\begin{figure}
   \centering
   \includegraphics[width=4in]{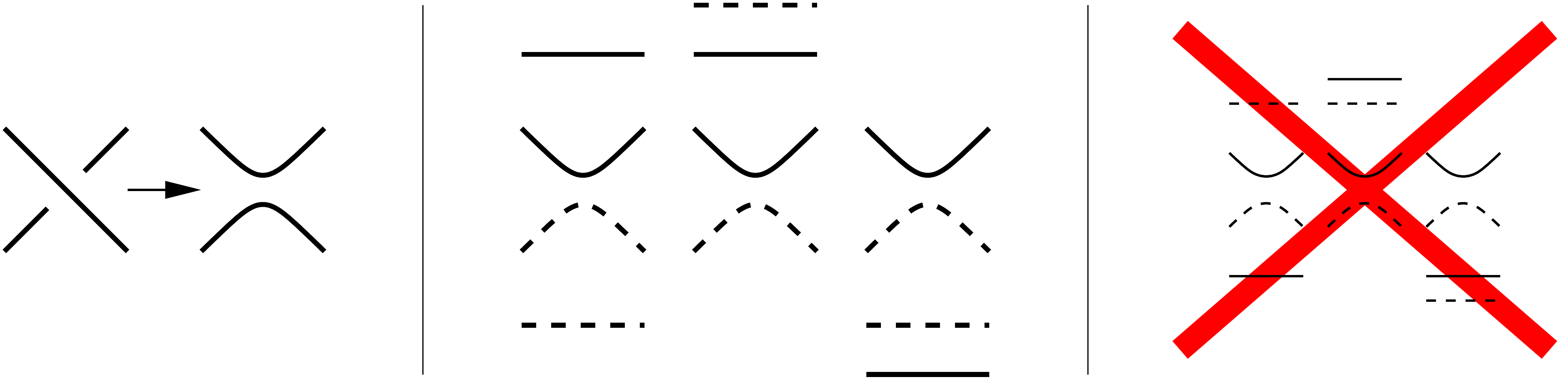}
   \caption{Allowed and disallowed configurations for switches of rulings}
   \label{fig:switch}
\end{figure}

\begin{Def}\label{def:rul} A \emph{ruling} is a partial smoothing
of a front diagram $f$ where certain crossings, called \emph{switches}, are replaced by a pair of arcs as in Figure \ref{fig:switch} so that the diagram becomes a union of standard unknot diagrams, called \emph{eyes}. (An eye is a pair of arcs connecting the same pair of left and right cusps that contain no other cusps and that otherwise do not meet, not even at switches.) We assume the so-called \emph{normality condition}: in the vertical ($x=\text{const.}$) slice of the diagram through each switch, the two eyes that meet at the switch fit one of the three configurations in the middle of Figure \ref{fig:switch}. 

The notion above is also known as an \emph{ungraded ruling}. 
If we assume that all switches are of even index, we get \emph{$2$--graded rulings}.\footnote{In the multi-component case, orientation of the diagram and our `even-right' convention for the Maslov potential is important.} The set of these objects will be denoted by $\Gamma_2(f)$. Finally, \emph{$\Z$--graded rulings} 
are those in which each switch has index $0$.
\end{Def}

\begin{Def}
The eyes bound disks and with the addition of twisted bands at the switches, we
may consider the ruling as a surface.\footnote{Compare with Seifert's construction and with the graph $K_D$ of \cite[\S6]{chp}.} This may be referred to as the \emph{surface associated to the ruling}, but usually we will simply identify the ruling and the surface. The \emph{genus} of a ruling is the genus of its associated surface. 
\end{Def}

This surface does not come
with an embedding into $\R^3$ though. (In \cite{cob}, we will associate to it an immersed exact Lagrangian in $B^4$.) In fact it follows from the
results of this paper that such an embedding often does not exist. But first, let us make several simple observations.

\begin{all}\label{iranyit}
The surface associated to a $2$--graded ruling is orientable. If an ungraded ruling of a Legendrian link is topologically an orientable surface, then the link has an orientation so that the ruling is $2$--graded. 
\end{all}

\begin{proof}
Consider the Maslov potential on the front diagram.
The entire top arc of each eye in a $2$--graded ruling is either on odd or on even
potential and the entire bottom arc has the opposite parity.
Switches connect odd arcs to odd arcs and even arcs to even arcs.
Hence if we give `odd top' eyes one orientation of the projection
plane and to `even top' eyes the opposite orientation, then (by an examination of Figure \ref{fig:switch}) this
will extend continuously over the bands.

For the other claim, simply induce an orientation of the Legendrian from that of its ruling. Doing so, all switches will be positive crossings.
\end{proof}

Now the following is obvious: 

\begin{kov}
For front diagrams of (single-component) knots, a ruling is $2$--graded if and only if its associated surface is orientable. 
\end{kov}


\begin{megj} 
Ungraded, that is not-necessarily-orientable
versions of Theorem \ref{thm:canon} and Corollary \ref{kov} fail even for alternating knots. The knot $5_1$ (the positive $(5,2)$ torus knot) bounds an embedded M\"obius band, thus has unorientable genus $1$, and it also has a ruling of oriented genus $2$ (i.e., unoriented genus $4$). 
\end{megj}


\begin{all}
If a ruling is topologically a disk, then it is $\Z$--graded. In particular its boundary Legendrian has $r=0$.
\end{all}


\begin{proof}
The eyes of such a ruling are connected by the switches in a tree-like fashion. Consider a leaf of this tree and notice that because of the way the Maslov potential is defined, its single switch must be of index $0$. Then remove this eye from the diagram and induct on the number of eyes.
\end{proof}

\begin{megj}
In fact, from the existence of a $2$--graded ruling it
already follows that $r=0$ \cite{josh}. On the other hand, there exist Legendrians with $r=0$ and ungraded rulings, hence with maximal $tb$ in their smooth type \cite{rulpoly}, yet without a single $2$--graded ruling. One such example is the untwisted double of the negative trefoil shown in Figure \ref{untwdblnegtref}, produced using the method of \cite{meginten} on a $+$adequate diagram of the knot. From Rutherford's results, to be reviewed in the next section, it follows that no Legendrian representative of this knot can have a $2$--graded ruling.

\begin{figure}\centering\includegraphics[width=3in,height=1.5in]{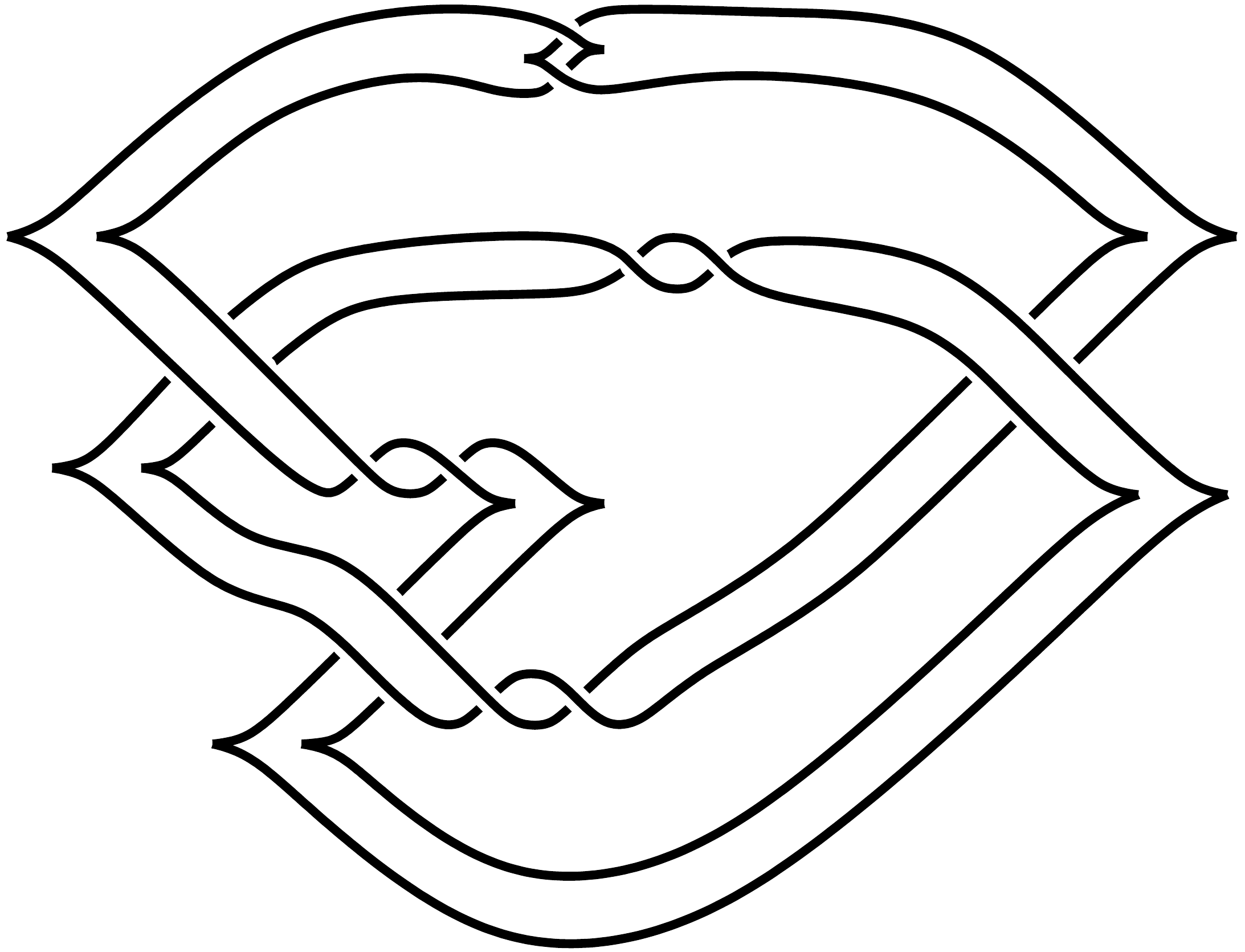}\caption{The untwisted double of the negative trefoil knot}\label{untwdblnegtref}\end{figure}

\end{megj}

\section{\small Proof of Theorem \ref{thm:canon}}\label{sec:dan}

We may classify rulings according to their genera. In fact the quantity \[\theta(\eta)=\text{number of eyes}-\text{number of
switches},\] introduced in \cite{chp}, is just the Euler characteristic of the surface associated to the ruling $\eta$. 


The counts of ungraded, $2$--graded,
and $\Z$--graded rulings with a given genus (or Euler characteristic) are
Legendrian isotopy invariants \cite{chp}. The sequence of these counts for all
genera 
is called
the \emph{complete ruling invariant}.  
In the $\Z$--graded case, the complete ruling invariant is
very effective in distinguishing Legendrian knots with the same
classical invariants. By contrast, the ungraded and $2$--graded
complete ruling invariants are determined by the smooth type and the
Thurston--Bennequin number. This is true by the two main theorems of
\cite{rulpoly}. In particular, in the $2$--graded case Rutherford proves the following. 

\begin{tetel}[\cite{rulpoly}]\label{thm:dan}
Let $P_L(v,z)$ denote the Homfly polynomial of the
Legendrian knot $L$, and let $Q_L(z)$ be the coefficient of
$v^{tb(L)+1}$ in $P_L(v,z)$. Then
\[Q_L(z) = \sum_{\eta\in\Gamma_2(L)} z^{2g(\eta)}.\]
\end{tetel}

\begin{megj}\label{max}
On the other hand, the smallest $v$--degree in the Homfly polynomial, which we will denote by $e$, is a well-known strict upper bound on $tb$ \cite{fw,morton,ft}. By the above, the existence of a $2$--graded ruling implies that this bound is sharp (that is, $tb+1=e$), i.e.\ that the Legendrian in question has maximal Thurston--Bennequin number within its smooth isotopy class. Conversely, $tb+1=e$ implies that our Legendrian does have $2$--graded rulings.
\end{megj}

In Theorem \ref{thm:dan}, our replacement of Rutherford's \[j(\eta)=\text{number of
switches}-\text{number of left cusps}+1\] with $2g(\eta)$ is valid 
by Proposition \ref{iranyit}. Note that we also replaced the variable $a$ in his formulation of the Homfly polynomial by $v^{-1}$.

\begin{Def}
For a knot type $K$, we define the \emph{ruling genus of $K$}, denoted by $\rho(K)$, to be the supremum of the genera of all $2$--graded rulings of all Legendrian representatives of $K$.
\end{Def}

This is either a finite value (and in that case a maximum) or $-\infty$. Note that by Theorem \ref{thm:dan} and Remark \ref{max}, any front diagram with smooth type $K$ and maximal Thurston--Bennequin number may be used to determine $\rho(K)$. In particular, if we know that $\rho\ne-\infty$, then we may directly read $\rho$ off of the Homfly polynomial. On the other hand, deciding whether $\rho(K)=-\infty$ is essentially equivalent to knowing the maximum Thurston--Bennequin number of $K$. At the time of this writing, that is still a hard problem in general.

With this new terminology, Corollary \ref{kov} simply asserts that
\[\rho(K)\le g(K)\text{ for all homogeneous }K,\] whereas by Theorem \ref{thm:canon}, \[\rho(K)\le \tilde g(K)\text{ for all }K.\]

\begin{proof}[Proof of Theorem \ref{thm:canon}]
From Rutherford's theorem above, it is obvious that the ruling genus $\rho$ is at most half of the $z$--degree of the Homfly polynomial, which in turn is a well known lower bound \cite[Theorem 2]{morton} for the canonical genus of the knot. 
\end{proof}

We close this section with one last simple observation, which essentially follows from \cite[Proposition 11.2]{chp}.

\begin{all}\label{pro:connsum}
The ruling genus is additive for connected sums: \[\rho(K_1\#K_2)=\rho(K_1)+\rho(K_2).\]
\end{all}

\begin{proof}
\begin{figure}
   \centering
   \includegraphics[width=4in]{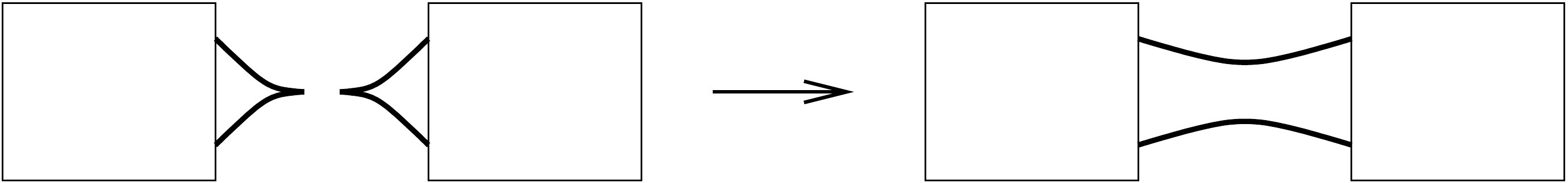}
   \caption{Connected sum operation}
   \label{fig:connsum}
\end{figure}
Choose front diagrams for $K_1$ and $K_2$ with maximal Thurston--Benne\-quin number. One may form their connected sum as in Figure \ref{fig:connsum}. Assume without loss of generality that the Maslov potentials near the two merging cusps match. Then it is plain to see that with any of the three notions of Definition \ref{def:rul}, the set of rulings for the new diagram is the Cartesian product of the two old sets. The claim follows at once.
\end{proof}

\section{\small Concluding remarks}\label{sec:jim}


One finds several sufficient conditions for $\rho\le g$ in terms of polynomial knot invariants.\footnote{For a long time, I tried to prove that $\rho\le g$ for all knots. Whitehead doubles, having genus $1$, are natural candidates for a counterexample; see Example \ref{ellenpelda}. However I initially dismissed them, having misunderstood a claim in \cite{fuchs}.} These are listed in the two propositions below. As before, let $e$ denote the minimum $v$--degree of the Homfly polynomial $P_K(v,z)$ of $K$. Versions of \eqref{negative} and \eqref{subset} on the following list were also published in \cite{private}.

\begin{all}\label{pro:noruling}
For a smooth knot type $K$, any of the following conditions implies that $K$ doesn't have Legendrian representatives with $2$--graded rulings, i.e.\ that $\rho(K)=-\infty$. 
\begin{enumerate}[(i)]
\item\label{lenny} ``Khovanov beats Homfly'': $e\ge2+\min\{\,k\,\big|\,\bigoplus_{i-j=k}HKh^{i,j}(K)\ne0\,\}$, where $HKh^{i,j}(K)$ is the Khovanov homology group of $K$ in bigrading $(i,j)$.
\item\label{kauff} ``Kauffman beats Homfly'': The minimum $v$--degree of the Kauffman polynomial of $K$ is less than $e$.
\item\label{negative} ``Negative counts'':  there is an integer $i$ so that the coefficient $p_{e,i}$ of $v^e z^i$ in the Homfly polynomial is negative.
\item\label{subset} ``Subset failure'': The minimum $v$-degree in the Kauffman polynomial is $e$ or more, but there is an integer $i$ so that if $f_{e,i}$ is the coefficient of $v^ez^i$ in the Dubrovnik version of the Kauffman polynomial, then $0\le p_{e,i}\le f_{e,i}$ fails to hold.
\end{enumerate}
\end{all}

\begin{proof}
Recall from Remark \ref{max} that $2$--graded rulings exist if and only if $tb+1=e$ for some Legendrian representing $K$, so this is what we have to prevent. 

In the first two cases, we separate $tb+1$ and $e$ by inequalities. Claim \eqref{lenny} is a direct consequence of the so-called Khovanov bound of the Thurston--Bennequin number \cite{khov}. Claim \eqref{kauff} follows in the same straightforward way from the so-called Kauffman bound \cite{congr,ft}. 

In the other two cases, we assume $tb+1=e$ and find a contradiction with Rutherford's results. Claim \eqref{negative} is obvious from Theorem \ref{thm:dan}. To prove the last statement, one needs to compare Theorem \ref{thm:dan} with Rutherford's other main theorem in \cite{rulpoly}, which says that the coefficients $f_{tb+1,i}$ count ungraded rulings of some given ($\theta=1-i$, to be exact) Euler characteristic. Now the contradiction is obvious from the fact that $2$--graded rulings are also ungraded.
\end{proof}

Let $M$ denote the maximum $z$--degree that appears in the Homfly polynomial of $K$ in a monomial that also contains $v^e$. With this, either $\rho=M/2$ or $\rho=-\infty$. (It is difficult to tell which one is true in the case of an arbitrary knot, even if none of the conditions listed in Proposition \ref{pro:noruling} holds. At present, the best one can do is to try and construct a Legendrian representative with $tb+1=e$. If such a front is found, then we know that $\rho=M/2$.)


\begin{all}
For a smooth knot type $K$, any of the following conditions implies that $\rho(K)\le g(K)$.
\begin{enumerate}[(a)]
\item ``Bennequin test'': For the exponents defined above, we have $M\le e$.
\item ``Conway test'': The degree of the Conway polynomial\hspace{3pt}\footnote{Recall that the Homfly polynomial reduces to the Conway polynomial by the substitution $\triangledown_K(z)=P_K(1,z)$.} $\triangledown_K(z)$ is $M$ or more.
\end{enumerate}
\end{all}

\begin{proof}
There is nothing to prove if $\rho=-\infty$, so we will assume that $\rho=M/2$. That implies $tb+1=e$ for those Legendrian representatives that have $2$--graded rulings. Now, the Bennequin test is so called because it is obvious from \eqref{benn}:
\[\rho=\frac M2\le \frac e2=\frac{tb+1}2\le g.\]

To prove the Conway test, recall that half the degree of the Conway polynomial is a lower bound for the genus, so in this case $\rho\le M/2\le g$.
\end{proof} 

Using the conditions listed above, Jim Hoste kindly conducted a computer search of the 
table of knots for potential examples with $\rho>g$. Knots with $M=2$ (i.e., $\rho=1$ or $-\infty$) were also ruled out. Among prime knots up to $13$ crossings, only two candidates were found: $13n_{1426}$ and $13n_{1456}$. It turns out that these knots do have $2$--graded rulings; in fact $\rho=2$ for both, with unique genus $2$ rulings. In both cases, I was able to add a $1$--handle to that ruling, realize the result as an embedded surface, and then find a Gabai disk decomposition \cite{gabai}. Thus both knots have $g=3$. However at crossing number $14$, the phenomenon of Example \ref{ellenpelda} occurs.

\end{document}